\documentclass{amsart}

\newtheorem{theorem}{Theorem}
\newcommand{\bt}{\begin{theorem}}
\newcommand{\et}{\end{theorem}}

\newtheorem*{theoremNN}{Theorem}
\newcommand{\btNN}{\begin{theoremNN}}
\newcommand{\etNN}{\end{theoremNN}}

\newtheorem{lemma}{Lemma}
\newcommand{\bl}{\begin{lemma}}
\newcommand{\el}{\end{lemma}}

\newtheorem{corollary}{Corollary}
\newcommand{\bc}{\begin{corollary}}
\newcommand{\ec}{\end{corollary}}

\newtheorem{definition}{Definition}
\newcommand{\bdf}{\begin{definition}}
\newcommand{\edf}{\end{definition}}

\newtheorem{conjecture}{Conjecture}
\newcommand{\bconj}{\begin{conjecture}}
\newcommand{\econj}{\end{conjecture}}

\newtheorem*{conjectureNN}{Conjecture}
\newcommand{\bconjNN}{\begin{conjectureNN}}
\newcommand{\econjNN}{\end{conjectureNN}}

\newtheorem{example}{Example}
\newcommand{\bex}{\begin{example}}
\newcommand{\eex}{\end{example}}

\newtheorem{problem}{Problem}
\newcommand{\bprob}{\begin{problem}}
\newcommand{\eprob}{\end{problem}}

\newtheorem*{problemNN}{Problem}
\newcommand{\bprobNN}{\begin{problemNN}}
\newcommand{\eprobNN}{\end{problemNN}}

\newtheorem{oproblem}{Open Problem}
\newcommand{\boprob}{\begin{oproblem}}
\newcommand{\eoprob}{\end{oproblem}}

\newtheorem*{oproblemNN}{Open Problem}
\newcommand{\boprobNN}{\begin{oproblemNN}}
\newcommand{\eoprobNN}{\end{oproblemNN}}

\newcommand{\beq}{\begin{equation}}
\newcommand{\eeq}{\end{equation}}
\newcommand{\benum}{\begin{enumerate}}
\newcommand{\eenum}{\end{enumerate}}

\newcommand{\Q}{\ensuremath{ \mathbf{Q} }}

\newcommand{\R}{\ensuremath{\mathbf R}}

\newcommand{\mbc}{\ensuremath{\mathbf c}}

\newcommand{\mbr}{\ensuremath{\mathbf r}}

\newcommand{\mcc}{\ensuremath{ \mathcal C}}

\newcommand{\mcr}{\ensuremath{ \mathcal R}}

\newcommand{\bq}{\begin{eqnarray*}}
\newcommand{\eq}{\end{eqnarray*}}
\newcommand{\be}{\begin{eqnarray}}
\newcommand{\ee}{\end{eqnarray}}
\newcommand{\ba}{\begin{array}}
\newcommand{\ea}{\end{array}}
\newcommand{\bfr}{\begin{flushright}}
\newcommand{\efr}{\end{flushright}}

\newcommand{\bmat}{\left(\begin{matrix}}
\newcommand{\emat}{\end{matrix}\right)}
\newcommand{\bsmallmat}{\left(\begin{smallmatrix}}
\newcommand{\esmallmat}{\end{smallmatrix}\right)}

\DeclareMathOperator{\colsum}{\text{colsum}}

\DeclareMathOperator{\diag}{\text{diag}}

\DeclareMathOperator{\qand}{\quad\text{and}\quad}
\DeclareMathOperator{\qqand}{\qquad\text{and}\qquad}

\DeclareMathOperator{\rowsum}{\text{rowsum}}

\usepackage{amssymb,latexsym}

\title[Explicit matrix limits]{Matrix scaling and explicit doubly stochastic limits}
\author{Melvyn B. Nathanson}
\address{Department of Mathematics\\Lehman College (CUNY)\\Bronx, NY 10468}
\email{melvyn.nathanson@lehman.cuny.edu}



\thanks{Supported in part by a grant from the PSC-CUNY Research Award Program.}

\date{\today}

\begin{document}
\maketitle

\begin{abstract}
The process of alternately row scaling and column scaling a positive $n \times n$ matrix $A$ 
converges to a doubly stochastic positive $n \times n$ matrix $S(A)$, often called the 
\emph{Sinkhorn limit} of $A$.  The main result in this paper is the computation of exact formulae 
for the Sinkhorn limits of certain symmetric positive $3\times 3$ matrices.  
\end{abstract}


\section{Doubly stochastic matrices and scaling}
Let $A = (a_{i,j})$ be an $n \times n$ matrix.
For $i \in \{1,\ldots, n\}$, the $i$th \emph{row sum}\index{row sum} of $A$ is 
\[
\rowsum_i(A) = \sum_{j=1}^n a_{i,j}.
\]
For $j \in \{1,\ldots, n\}$, the $j$th \emph{column sum}\index{column sum} of $A$ is 
\[
\colsum_j(A) = \sum_{i =1}^n a_{i,j}. 
\]
The matrix  $ A = (a_{i,j})$ is \emph{positive} if $a_{i,j}>0$ for all $i$ and $j$, 
and \emph{nonnegative} if $a_{i,j} \geq 0$ for all $i$ and $j$. 
The matrix $A = (a_{i,j})$ is \emph{row stochastic}\index{row stochastic} 
if $A$ is nonnegative and  $\rowsum_i(A) = 1$ for all $i \in \{1,\ldots, n\}$.  
The  matrix $A$ is \emph{column stochastic}\index{column stochastic} 
if $A$ is nonnegative and  $\colsum_j(A) = 1$ for all $j \in \{1,\ldots, n\}$.  
The matrix $A$ is \emph{doubly stochastic}\index{doubly stochastic} 
if it is both row and column stochastic.  


Let  $\diag(x_1,\ldots, x_n)$ denote the $n \times n$ diagonal matrix 
whose $(i,i)$th coordinate is $x_i$ for  all $i \in \{1,2,\ldots, n\}$.  
The matrix $\diag(x_1,x_2,\ldots, x_n)$ 
is \emph{positive  diagonal} if $x_i > 0$ for all $i$.

Let $A = (a_{i,j})$ be an $n \times n$ matrix.  
The process of multiplying the rows of $A$ by scalars, 
or, equivalently, multiplying $A$ on the left by a diagonal matrix $X$, 
is called \emph{row-scaling}, 
and $X$ is called a \emph{row-scaling matrix}.  

The process of multiplying the columns of $A$ by scalars, 
or, equivalently, multiplying $A$ on the right by a diagonal matrix $Y$,  
is called \emph{column-scaling}, 
and $Y$ is called a \emph{column-scaling matrix}.  

If $X = \diag(x_1,x_2,\ldots, x_n)$ and $Y = \diag(y_1,y_2,\ldots, y_n)$, then 
\[
XAY = 
\bmat
x_1a_{1,1} y_1 & x_1a_{1,2} y_2 & x_1a_{1,3} y_3 &  \cdots & x_1a_{1,n} y_n \\ 
x_2a_{2,1} y_1 & x_2a_{2,2} y_2 &x_2a_{2,3} y_3 &  \cdots & x_2a_{2,n} y_n \\ 
\vdots &&&& \vdots \\
x_n a_{n,1} y_1 & x_n a_{n,2} y_2 &x_n a_{n,3} y_3 &  \cdots & x_n a_{n,n} y_n 
\emat.
\]

Let $A = (a_{i,j})$ be an $n \times n$ matrix with positive row sums, 
that is, $\rowsum_i(A) > 0$ for all $i \in \{1,\ldots, n\}$.  
Let 
\beq                     \label{Sinkhorn:rowsum}
X(A) = \diag\left( \frac{1}{\rowsum_1(A)}, \ldots,  \frac{1}{\rowsum_n(A)} \right)
\eeq 
and let 
\[
\mcr(A) = X(A) A.
\]
We have 
\[
\mcr(A)_{i,j} = \frac{a_{i,j}}{\rowsum_i(A)}
\]
and so
\[
\rowsum_i(\mcr(A)) 
= \sum_{j=1}^n \mcr(A)_{i,j} 
= \sum_{j=1}^n\frac{a_{i,j}}{\rowsum_i(A)} 
= \frac{\rowsum_i(A)}{\rowsum_i(A)} = 1
\]
for all $i \in \{1,\ldots, n \}$.  
Therefore, $\mcr(A)$ is a row stochastic matrix.

Similarly, if $A = (a_{i,j})$ is an $n \times n$ matrix with positive column sums and if 
\beq                     \label{Sinkhorn:colsum}
 Y(A) = \diag\left( \frac{1}{\colsum_1(A)}, \ldots,  \frac{1}{\colsum_n(A)} \right)
 \eeq
and 
\[
\mcc(A) = A Y(A),
\]
then 
\[
\mcc(A)_{i,j} = \frac{a_{i,j}}{\colsum_j (A)}
\]
and 
\[
\colsum_j(\mcc(A)) 
= \sum_{i=1}^n \mcc(A)_{i,j} = \sum_{i=1}^n \frac{a_{i,j}}{\colsum_j(A)} 
= \frac{\colsum_j(A)}{\colsum_j(A)} = 1
\]
for all $j \in \{1,\ldots, n\}$.  
Therefore, $\mcc(A)$ is a column stochastic matrix.

The following two theorems were stated by Sinkhorn~\cite{sink64}, 
and subsequently proved by 
Brualdi, Parter, and Schneider~\cite{brua-part-schn66},
Djokovi\'{c}~\cite{djok70}, 
Knopp-Sinkhorn~\cite{sink-knop67}, 
Menon~\cite{meno67},  Letac~\cite{leta74}, and Tverberg~\cite{tver76}.

\bt                 \label{Sinkhorn:theorem:Sinkhorn-1}
Let $ A = (a_{i,j})$ be a positive $n \times n$ matrix.
\benum
\item[(i)]
There exist positive diagonal $n\times n$ matrices $ X$ and $ Y$ 
such that $ X A  Y$ is doubly stochastic.
\item[(ii)]
If $ X$, $ X'$, $ Y$, and $ Y'$ are  positive diagonal $n\times n$ matrices such that 
both  $ X A  Y$ and $ X' A  Y'$ are doubly stochastic, then 
$ X A  Y =  X' A  Y'$ and there exists $\lambda > 0$ 
such that $ X' = \lambda  X$ and $ Y' = \lambda^{-1}  Y$.  
\item[(iii)]
Let A\ be a positive symmetric $n \times n$ matrix.  There exists a unique positive diagonal matrix X\ 
such that $X A X$ is doubly stochastic. 
\eenum
\et

The unique doubly stochastic matrix $XAY$ in Theorem~\ref{Sinkhorn:theorem:Sinkhorn-1} 
is called the \emph{Sinkhorn limit} of A, and denoted $S(A)$.

\bt                             \label{Sinkhorn:theorem:Sinkhorn-2}
Let A\ be a positive $n\times n$ matrix, and let $S(A)$ be the Sinkhorn limit of $A$.
Construct  sequences of positive matrices  $(A_{\ell})_{\ell=0}^{\infty}$ 
and $(A'_{\ell})_{\ell=0}^{\infty}$ 
and sequences of positive diagonal matrices $(X_{\ell})_{\ell=0}^{\infty}$ 
and $(Y_{\ell})_{\ell=0}^{\infty}$ 
as follows:
Let 
\[
A_0 = A.
\]
Given the matrix $A_{\ell}$, let  
\beq                                               \label{Sinkhorn:X}
X_{\ell} = X(A_{\ell}) 
\eeq 
be the row-scaling matrix of $A_{\ell}$ defined by~\eqref{Sinkhorn:rowsum}. 
The matrix 
\[
A'_{\ell} = \mcr(A_{\ell}) =  X_{\ell} A_{\ell}.
\]
is row stochastic.  
Let 
\beq                                               \label{Sinkhorn:Y}
Y_{\ell} = Y(A'_{\ell}) 
\eeq
be the column-scaling matrix of $A'_{\ell}$ defined by~\eqref{Sinkhorn:colsum}, and let 
\[
A_{\ell+1} =  \mcc(A'_{\ell}) =   A'_{\ell}Y_{\ell}.
\]
The matrix $A_{\ell +1}$ is column stochastic.  

The Sinkhorn limit is obtained by alternately row-scaling and column-scaling:
\[
S(A) = \lim_{\ell\rightarrow \infty} A_{\ell} = \lim_{\ell\rightarrow \infty} A'_{\ell}.
\] 
\et

It is an open problem to compute explicitly the Sinkhorn limit of a positive $n\times n$ matrix.  
This is known  for $2\times 2$ matrices (Nathanson~\cite{nath2019-184}).  
The goal of this paper is the explicit computation of Sinkhorn limits 
for certain $3 \times 3$ matrices.


\section{Sinkhorn limits of $3 \times 3$  symmetric matrices  and 
their doubly stochastic shapes}

Let $A$ and $B$ be positive  $n\times n$ matrices.
We write $A \sim B$ if there exist $n\times n$  permutation matrices $P$ and $Q$ 
and $\lambda > 0$ such that 
\[
B = \lambda PAQ.
\]
This is an equivalence relation.  Moreover,    
$A\sim B$ implies 
\beq                                 \label{Sinkhorn:PSQ}
S(B) = \lambda P S(A) Q.
\eeq
Thus, it suffices to determine the Sinkhorn limit of only one matrix in an equivalence class.

We shall compute the Sinkhorn limit of every symmetric positive $3\times 3$ matrix
whose set of coordinates consists of two distinct real numbers.
  
Let $A$ be such a matrix with coordinates $M$ and $N$ with $M \neq N$. 
There are 9 coordinate positions in the matrix, and so exactly one 
of the numbers $M$ and $N$ occurs at least five times.  
Suppose that the coordinate $M$ occurs five or more times.  Let $\lambda = 1/M$ 
and $K = N/M$.
The matrix $\lambda A$ has two distinct positive  coordinates $1$ and $K$, 
and $K$ occurs at most four times.  
There are seven equivalence classes of such matrices with respect to permutations 
and dilations.
The main result of this paper is the calculation of the 
Sinkhorn limits of these matrices.    

\bt
Let $K >0$ and $K \neq 1$.  
The matrices $A_1,\ldots, A_7$ below are a complete set 
of representatives of the seven equivalence classes 
of symmetric $3 \times 3$ matrices with coordinates 1 and $K$.  
The matrix $S(A_i)$ gives the shape of the Sinkhorn limit of $A_i$ for $i = 1,\ldots, 7$. 
The coordinates of the Sinkhorn limits as explicit functions of 1 and $K$ are computed  
in Sections~\ref{Sinkhorn:section:A_1}--\ref{Sinkhorn:section:A_7}.
\benum
\item
\[
A_1 = \bmat
K & 1 & 1 \\
1 & K & 1 \\
1 & 1 & K
\emat
\qquad
S(A_1) = \bmat
a & b & b \\
b & a & b \\
b & b & a  
\emat
\]

\item
\[
A_2 = \bmat
K & 1 & 1 \\
1 & 1 & 1 \\
1 & 1 & 1
\emat
\qquad
S(A_2) = 
\bmat
a & b & b \\
b & c & c \\
b & c & c 
\emat
\]

\item
\[
A_3 = \bmat 
1 & 1 & 1 \\
1 & K & K  \\
1 & K & K  
\emat
\qquad
S(A_3) = \bmat
a & b & b \\
b & c & c \\
b & c & c
\emat
\]

\item
\[
A_4 = \bmat
1 & K & K \\
K & 1 & 1 \\
K & 1 & 1
\emat
\qquad
S(A_4) = \bmat
a & b & b \\
b & c & c \\
b & c & c
\emat
\]

\item
\[
A_5 = \bmat
K & 1 & 1 \\
1 & K & 1 \\
1 & 1 & 1
\emat
\qquad
S(A_5) = \bmat
a & b & c \\
b & a & c \\
c & c & d 
\emat
\]

\item
\[
A_6 = \bmat
K & K & 1 \\
K & 1 & 1 \\
1 & 1 & 1
\emat
\qquad
S(A_6) = \bmat
a & b &c \\
b & c & a \\
c & a & b 
\emat
\]

\item
\[
A_7 = \bmat
K  & K & 1 \\
K & 1 & 1  \\
1 & 1 & K
\emat
\qquad
S(A_7) = \bmat
a & b & c \\
b & d & e \\
c & e & f
\emat
\]

\eenum

\et


\section{The $MBN$ matrix}
Let $k$, $\ell$, and $n$ be positive integers such that 
\[
n = k+\ell. 
\]
Let $M$, $B$, and $N$ be positive real numbers.  
Consider the $n \times n$ symmetric matrix 
\beq            \label{Sinkhorn:MBN}
A = \bmat
M & M & \cdots & M & B & B & \cdots & B \\
M & M & \cdots & M & B & B & \cdots & B \\
\vdots &&& \vdots  & \vdots  &&& \vdots \\ 
M & M & \cdots & M & B & B & \cdots & B \\
 B & B & \cdots & B & N & N & \cdots & N \\
  B & B & \cdots & B & N & N & \cdots & N \\
\vdots &&& \vdots  & \vdots  &&& \vdots \\ 
 B & B & \cdots & B & N & N & \cdots & N \\
\emat
\eeq
in which the first $k$  rows are equal to 
\[
(\underbrace{M,M,\ldots, M}_{k},\underbrace{ B, B, \ldots, B}_{\ell})
\]
and the last $\ell$ rows  are equal to 
\[
(\underbrace{B,B,\ldots, B}_{k}, \underbrace{ N,N, \ldots, N}_{\ell}).
\]
Let $X  = \diag(x_1,x_2,x_3, \ldots, x_n ) $  be the unique positive $n \times n$  diagonal matrix 
such that the alternate scaling limit  $S(A) = X  A X $ is doubly stochastic. 
Thus, the matrix 
\[
S(A) 
 = \bmat
M x_1^2  & Mx_1x_2 & \cdots & M x_1x_k & B x_1x_{k+1} & B x_1x_{k+2} & \cdots & B x_1x_{n}\\
M x_2x_1  & Mx_2^2 & \cdots & M x_2x_k & B x_2x_{k+1} & B x_2x_{k+2} & \cdots & B x_2x_{n}\\
\vdots &&& \vdots  & \vdots  &&& \vdots \\ 
M x_k x_1  & Mx_kx_2 & \cdots & M x_k^2 & B x_kx_{k+1} & B x_kx_{k+2} & \cdots & B x_kx_{n}\\
B x_{k+1} x_1  & B x_{k+1} x_2 & \cdots & B x_{k+1}x_k & N x_{k+1}^2& N x_{k+1}x_{k+2} & \cdots & N x_{k+1}x_{n}\\
B x_{k+2} x_1  & B x_{k+2} x_2 & \cdots & B x_{k+2}x_k & N x_{k+2}x_{k+1}& N x_{k+2}^2 & \cdots & N x_{k+2}x_{n}\\
\vdots &&& \vdots  & \vdots  &&& \vdots \\ 
B x_n x_1  & B x_n x_2 & \cdots & B x_n x_k & N  x_nx_{k+1}& N x_nx_{k+2} & \cdots & N x_{n}^2
\emat
\]
satisfies 
\[
x_i\left( M \sum_{j=1}^k x_j + B\sum_{j=k+1}^n x_j \right) = 1 
\qquad \text{for $i = 1,2,\ldots k$} 
\]
and
\[
x_i\left( B \sum_{j=1}^k x_j + N\sum_{j=k+1}^n x_j \right) = 1 
\qquad \text{for $i = k+1, k+2, ,\ldots k+\ell$.} 
\]
It follows that $x_i = x_1$ for  $i = 1,2,\ldots k$ and  $x_i = x_n$ for  $i = k+1, k+2,\ldots k + \ell$.  
Let $x_1 = x$ and $x_{k+1} = y$.  Define the diagonal matrix 
\[
X = \diag(\underbrace{x,x,\ldots, x}_{k}, \underbrace{ y,y, \ldots, y}_{\ell}).
\] 
We obtain 
\begin{align}                
S(A)     \label{Sinkhorn:MBN-limit}  
& =  \bmat
M x^2& M x^2& \cdots & Mx^2& B  xy & B xy & \cdots & B xy \\
M x^2& M x^2& \cdots & M x^2 & B xy & B xy& \cdots & B xy \\
\vdots &&& \vdots  & \vdots  &&& \vdots \\ 
M x^2 & Mx^2& \cdots & M x^2 & B  xy & B  xy & \cdots & B xy \\
B  xy & B  xy& \cdots & B  xy & N y^2  & N  y^2 & \cdots & N  y^2 \\
B  xy  & B  xy& \cdots & B  xy & N y^2 & N  y^2  & \cdots & N  y^2 \\
\vdots &&& \vdots  & \vdots  &&& \vdots \\ 
B  xy & B  xy & \cdots & B xy & N   y^2 & N  y^2  & \cdots & N y^2
\emat                           \\
& =
\bmat
a  & a & \cdots & a & b & b & \cdots & b \\
a  & a & \cdots & a & b & b & \cdots & b \\
\vdots &&& \vdots  & \vdots  &&& \vdots \\ 
a  & a & \cdots & a & b & b & \cdots & b \\
 b & b & \cdots & b & c & c & \cdots & c \\
 b & b & \cdots & b &  c & c & \cdots & c \\
\vdots &&& \vdots  & \vdots  &&& \vdots \\ 
 b & b & \cdots & b & c & c & \cdots & c              \nonumber
 \emat                       
\end {align}
where 
\begin{align}       \label{Sinkhorn:MBN-a}
a & = Mx^2                            \\ 
b & = Bxy  = \frac{1-ka}{\ell}               \label{Sinkhorn:MBN-b}     \\ 
c & = Ny^2  = \frac{1-kb}{\ell} =  \frac{\ell-k + k^2 a}{\ell^2}.          \label{Sinkhorn:MBN-c}    
\end{align}
Because $S(A)$ is row stochastic, we have 
\beq            \label{Sinkhorn:MBN-1}
x \left(kMx+ \ell B y \right) = 1 
\eeq
and
\beq            \label{Sinkhorn:MBN-2}
y \left( kBx + \ell N y \right) = 1. 
\eeq
Equation~\eqref{Sinkhorn:MBN-1} gives 
\[
y =\frac{1}{\ell B} \left( \frac{1}{x} -kMx \right).
\]
Inserting this into equation~\eqref{Sinkhorn:MBN-2} and rearranging gives  
\beq            \label{Sinkhorn:MBN-3}
k^2 M \left( MN - B^2 \right) x^4 - \left( 2k (MN - B^2) + nB^2 \right) x^2+ N = 0 
\eeq 

If $MN-B^2 = 0$, then 
\[
x^2 = \frac{N}{nB^2} = \frac{1}{nM} 
\]
and $M x^2 = a = b = c = 1/n$.  Thus, $S(A)$ is the $n \times n$ doubly stochastic matrix 
with every coordinate equal to $1/n$.

If $MN - B^2 \neq 0$, then~\eqref{Sinkhorn:MBN-3} is a quadratic equation in $x^2$. 
Let 
\[
L = \frac{MN}{B^2}.
\]
We obtain 
\begin {align*}
x^2 
& = \frac{ 2k (MN - B^2)  + nB^2 \pm B\sqrt{ 4k\ell (MN-B^2) + n^2B^2}}{ 2k^2 M (MN-B^2)} \\
& =  \frac{ 1}{kM} + \frac{ nB^2 \pm B\sqrt{ 4k\ell MN  + (k-\ell)^2 B^2}}{ 2k^2 M (MN-B^2)} \\
& =  \frac{ 1}{kM} + \frac{ n \pm \sqrt{ 4k\ell L  + (k-\ell)^2 }}{ 2k^2 M( L-1)} \\
& =  \frac{ 1}{kM} + \frac{ n \pm \sqrt{ n^2 + 4k\ell (L-1) }}{ 2k^2 M( L-1)}
\end {align*}
and
\begin {align*}
a  = M x^2 
& =  \frac{ 1}{k} + \frac{ n \pm \sqrt{ n^2 + 4k\ell (L-1) }}{ 2k^2( L-1)}.  
\end {align*}
Recall that $ka+\ell b = 1$ and so 
\[
0 < a< \frac{1}{k}.
\]
If $MN > B^2$, then $L>1$ and 
\[
 \sqrt{ n^2 + 4k\ell (L-1) } > n >  0.
\]
The inequality $a < 1/k$ implies that 
\beq                \label{Sinkhorn:MBN-aa}                        
a  =  \frac{ 1}{k} +  \frac{ n - \sqrt{ n^2 + 4k\ell (L-1) }}{ 2k^2( L-1)}.  
\eeq
If $MN < B^2$, then $0 < L < 1$ and 
\[
a  =  \frac{ 1}{k} - \frac{ n \pm \sqrt{ n^2 - 4k\ell (1-L) }}{ 2k^2( 1 - L)}.  
\]
Because 
\[
 \frac{ n+  \sqrt{ n^2 - 4k\ell (1-L) }}{ 2k^2( 1 - L)} > \frac{1}{k}
\]
the inequality $a > 0$ implies~\eqref {Sinkhorn:MBN-aa}.

We have proved the following.    

\bt         \label{Sinkhorn:theorem:MBN-1-1}
The Sinkhorn limit of the $MBN$ matrix~\eqref{Sinkhorn:MBN} is a 
doubly stochastic matrix $S(A)$ with shape~\eqref{Sinkhorn:MBN-limit}.
If $L = MN/B^2 = 1$, then $a = b = c = 1/n$.  
If $L \neq 1$, then equations~\eqref{Sinkhorn:MBN-aa},~\eqref{Sinkhorn:MBN-b},  and~\eqref{Sinkhorn:MBN-c} define the coordinates $a$,  $b$, and $c$.
The  matrix $S(A)$ depends only on the ratio $MN/B^2$.
\et

For example, the matrices 
\[
\bmat 
2 & 5 & 5 \\
5 & 3 & 3 \\
5 & 3 & 3
\emat,
\qquad
\bmat 
6 & 5 & 5 \\
5 & 1 & 1 \\
5 & 1 & 1
\emat,
\qquad 
\bmat 
6/25 & 1 & 1 \\
1 & 1 & 1 \\
1 & 1 & 1
\emat 
\]
have the same Sinkhorn limit with 
\begin{align*} 
a & = - \frac{37}{38} +  \frac{5\sqrt{73}}{38} = 0.1505\ldots \\ 
b & =   \frac{75}{76} - \frac{5\sqrt{73}}{76} = 0.4247\ldots \\ 
c & =  \frac{1}{152} +  \frac{5\sqrt{73}}{152}  = 0.2876\ldots.
\end{align*}

Let $\left( A^{(r)} \right)_{r=1}^{\infty}$ be a sequence of $MBN$ matrices 
such that $\lim_{r\rightarrow \infty} MN/B^2 = \infty$.
Let 
\[
S\left( A^{(r)} \right) = 
\bmat
a^{(r)}  & a^{(r)} & \cdots & a^{(r)} & b^{(r)} & b^{(r)} & \cdots & b^{(r)} \\
a^{(r)}  & a^{(r)} & \cdots & a^{(r)} & b^{(r)} & b^{(r)} & \cdots & b^{(r)} \\
\vdots &&& \vdots  & \vdots  &&& \vdots \\ 
a^{(r)}  & a^{(r)} & \cdots & a^{(r)} & b^{(r)} & b^{(r)} & \cdots & b^{(r)} \\
 b^{(r)} & b^{(r)} & \cdots & b^{(r)} & c^{(r)} & c^{(r)} & \cdots & c^{(r)} \\
 b^{(r)} & b^{(r)} & \cdots & b^{(r)} &  c^{(r)} & c^{(r)} & \cdots & c^{(r)} \\
\vdots &&& \vdots  & \vdots  &&& \vdots \\ 
 b^{(r)} & b^{(r)} & \cdots & b^{(r)} & c^{(r)} & c^{(r)} & \cdots & c^{(r)} \\
 \emat.
\]
We have 
\[
\lim_{r \rightarrow \infty} a^{(r)} = \frac{1}{k}, \qquad
\lim_{r \rightarrow \infty} b^{(r)}  = 0,   
\qquad 
\lim_{r \rightarrow \infty} c^{(r)} = \frac{1}{\ell}
\]
and
\[
\lim_{r\rightarrow \infty} S\left( A^{(r)} \right) = 
\bmat
1/k  & 1/k & \cdots & 1/k & 0 & 0 & \cdots & 0 \\
1/k  & 1/k & \cdots & 1/k & 0 & 0 & \cdots & 0 \\
\vdots &&& \vdots  & \vdots  &&& \vdots \\ 
1/k  & 1/k & \cdots & 1/k & 0 & 0 & \cdots & 0 \\
 0 & 0 & \cdots & 0 & 1/\ell & 1/\ell & \cdots & 1/\ell \\
 0 & 0 & \cdots & 0 &  1/\ell & 1/\ell & \cdots & 1/\ell \\
\vdots &&& \vdots  & \vdots  &&& \vdots \\ 
 0 & 0 & \cdots & 0 & 1/\ell & 1/\ell & \cdots & 1/\ell 
 \emat.
\]
Similarly, let $\left( A^{(r)} \right)_{r=1}^{\infty}$ be a sequence of $MBN$ matrices 
such that $\lim_{r\rightarrow \infty} MN/B^2 = 0$.  
It follows from~\eqref{Sinkhorn:MBN-a}  that 
\[
\lim_{r \rightarrow \infty} a^{(r)} 
=  \frac{ 1}{k} - \frac{ k+\ell  - |k-\ell | }{ 2k^2},
\]
If $k \leq \ell$, then 
\[
\lim_{r \rightarrow \infty} a^{(r)} = 0, 
\qquad 
\lim_{r \rightarrow \infty} b^{(r)} = \frac{1}{\ell}, 
\qquad
\lim_{r \rightarrow \infty} c^{(r)} =  \frac{\ell-k}{\ell^2}.  
\]
If $k > \ell$ , then 
\[
\lim_{r \rightarrow \infty} a^{(r)} =  \frac{k - \ell}{k^2},  
\qquad 
\lim_{r \rightarrow \infty} b^{(r)} =  \frac{1}{k}, 
\qquad
\lim_{r \rightarrow \infty} c^{(r)} = 0.
\]

\section{The matrix $A_1$}          \label{Sinkhorn:section:A_1} 
The matrix 
\[
A_1 = \bmat
K & 1 & 1 \\
1 & K & 1 \\
1 & 1 & K
\emat
\]
is the simplest.  Just one row scaling or one column scaling produces 
the doubly stochastic matrix 
\[
S(A_1) = \bmat
K/(K+2)  & 1/(K+2) & 1/(K+2) \\
1/(K+2) & K/(K+2) & 1/(K+2) \\
1/(K+2) & 1/(K+2) & K/(K+2)  
\emat
\]
We have  $S(A_1) = X A_1 X$, where 
\[
X = \diag( \sqrt{1/(K+2)}, \sqrt{1/(K+2)}, \sqrt{1/(K+2)} ).
\]
We have the asymptotic limits  
\[
\lim_{K\rightarrow \infty} S(A_1) = 
\bmat 
1 & 0 & 0 \\
0 & 1 & 0 \\
0 & 0 & 1 
\emat
\qand
\lim_{K\rightarrow 0} S(A_1) = 
\bmat 
0 & 1/2 & 1/2 \\
1/2 & 0 & 1/2 \\
1/2 & 1/2 & 0 
\emat.
\]

\section{The matrices $A_2$, $A_3$,  and $A_4$}             \label{Sinkhorn:section:A_2} 
These  are $MBN$ matrices.  
The matrix 
\[
A_2 = \bmat K & 1 & 1 \\ 1 & 1 & 1 \\ 1 & 1 & 1 \emat 
\] 
is an $MBN$ matrix with $k=1$, $\ell = 2$, $M=K$,  $B = N = 1$, and $L = K$.

The matrix 
\[
A_3 = \bmat 
1 & 1 & 1 \\
1 & K & K  \\
1 & K & K  
\emat
\]
is an $MBN$ matrix with $k=1$, $\ell = 2$, $M=B=1$,  $N = K$, and $L = K$.
Both matrices satisfy $L = MN/B^2 = K \neq 1$, and so they have the same Sinkhorn limit
\[
S(A_2) = S(A_3) = \bmat 
a & b & b \\
b & c & c \\
b & c & c
\emat
\]
with 
\begin{align}
a & =  \frac{2K+1 - \sqrt{8K  + 1}}{2(K-1)}     \label{Sinkhorn:A2-a}            \\
b& = \frac{ -3 + \sqrt{8K  + 1}}{4(K-1)}    \label{Sinkhorn:A2-b}            \\
c & = \frac{ 4K-1 - \sqrt{8K  + 1}}{8(K-1)}.   \label{Sinkhorn:A2-c}            
\end{align}

We have the asymptotic limits
\[
\lim_{K\rightarrow \infty} S(A_2) = 
\bmat 
1 & 0 & 0 \\
0 & 1/2 & 1/2 \\
0 & 1/2 & 1/2 
\emat
\qand
\lim_{K\rightarrow \infty} S(A_2) = 
\bmat 
0 & 1/2 & 1/2 \\
1/2 & 1/4 & 1/4 \\ 
1/2 & 1/4 & 1/4 
\emat.
\]

The matrix 
\[
A_4 = \bmat
1 & K & K \\
K & 1 & 1 \\
K & 1 & 1
\emat
\]
is an $MBN$ matrix with $k=1$, $\ell = 2$, $M=N=1$, and $B = K$.
We have $L = MN/B^2 = 1/K^2 \neq 0$, and the Sinkhorn limit 
\[
S(A_4) = 
\bmat 
a & b & b \\
b & c & c \\
b & c & c
\emat
\]
with 
\begin {align*}
a & =  \frac{-K^2 -2+ K \sqrt{K^2 + 8}}{2(K^2-1)} \\
b& =       \frac{3K^2 - K \sqrt{K^2 + 8}}{4(K^2-1)}      \\
c & =      \frac{K^2 - 4 + K \sqrt{K^2 + 8}}{8(K^2-1)}. 
\end {align*}

We have the asymptotic limits
\[
\lim_{K\rightarrow \infty} S(A_4) = 
\bmat 
0 & 1/2 & 1/2 \\
1/2 & 1/4 & 1/4 \\
1/2 & 1/4 & 1/4 
\emat
\qand
\lim_{K\rightarrow 0} S(A_4) = 
\bmat 
1 & 0 & 0 \\
0 & 1/2 & 1/2 \\
0 & 1/2 & 1/2 
\emat.
\]


\section{The matrix $A_5$}                    \label{Sinkhorn:section:A_5} 
The construction of the Sinkhorn limit of the $3 \times 3$ matrix
\[
A_5 = \bmat 
K & 1 & 1 \\
1 & K & 1 \\
1 & 1 & 1
\emat 
\]
requires only high school algebra.  
There exists a unique positive diagonal matrix $X = \diag(x,y,z)$  
such that $XA_5X$ is doubly stochastic and positive.  We have
\[
S(A_5) = XA_5X = \bmat 
Kx^2 & xy & xz \\
xy & Ky^2 & yz \\
xz & yz & z^2 
\emat
\]
and so 
\begin {align*}
Kx^2 + xy +  xz  & = 1 \\
xy + Ky^2 + yz  & = 1 \\
xz + yz + z^2   & = 1 
\end {align*}
We have 
\[
z = \frac{1 - Kx^2 - xy }{x} = \frac{1 - xy  - Ky^2 }{y}.  
\]
Rearranging, we obtain 
\beq                         \label{Sinkhorn:zC} 
(y-x)((K-1)xy+1) = 0.
\eeq
Note that $0 < xy < 1$. 
If $K > 1$, then  $(K-1)xy+1 > 1$.
 If $0 < K < 1$, then 
\[
0 < (1-K)xy < 1-K < 1
\]
and  $(K-1)xy+1 > 0$.
Therefore, $x=y$, and so 
\beq                         \label{Sinkhorn:zC} 
z = \frac{1 - (K+1)x^2}{x}
\eeq
\beq                         \label{Sinkhorn:zA} 
(K+1)x^2 + xz  = 1 
\eeq
\beq                      \label{Sinkhorn:zB} 
2xz + z^2  = 1.   
\eeq
We obtain   
\[
2\left( 1 - (K+1)x^2 \right)  +  \left(\frac{  1 - (K+1)x^2}{x} \right)^2 = 1.   
\]
Applying~\eqref{Sinkhorn:zC} and eliminating $xz$ from~\eqref{Sinkhorn:zA} and~\eqref{Sinkhorn:zB} gives 
\[
\left( \frac{1 - (K+1)x^2}{x} \right)^2 = z^2 = 2(K+1)x^2 - 1.
\]
Therefore, 
\[
(K^2-1)x^4 - (2K+1)x^2 + 1 = 0
\]
and so 
\[
x^2 = \frac{2K+1 \pm \sqrt{4K+5}}{2(K^2 -1)}.
\]
The inequality $Kx^2 < 1$ implies 
\[
x^2 = \frac{2K+1 - \sqrt{4K+5}}{2(K^2 -1)}
\]
and 
\[
z^2  = \frac{K+2 - \sqrt{4K+5}}{K-1}.
\]
Thus, the Sinkhorn limit has the shape 
\[
S(A_5) = \bmat
a & b & c \\
b & a & c \\
c & c & d
\emat
\]
where 
\begin {align*}
a & = Kx^2 = \frac{K(2K+1 - \sqrt{4K+5})}{2(K^2 -1)}  \\
b & = x^2  = \frac{2K+1 - \sqrt{4K+5}}{2(K^2 -1)}  \\
c & = xz = \frac{ \sqrt{ 2K+7 - 3\sqrt{4K+5}}}   {\sqrt{2}(K-1)}  \\
d & =   z^2  = \frac{K+2 - \sqrt{4K+5}}{K-1}.
\end {align*}

We have the asymptotic limits
\[
\lim_{K\rightarrow \infty} S(A_5) = \bmat 1 & 0 & 0 \\ 0 & 1 & 0 \\ 0 & 0 & 1 \emat 
\qand
\lim_{K\rightarrow 0} S(A_5) = \bmat 
0 & \frac{\sqrt{5}-1}{2} & \frac{3 - \sqrt{5}}{2}  \\  
\frac{\sqrt{5}-1}{2}   & 0 & \frac{3 - \sqrt{5}}{2} \\ 
\frac{3 - \sqrt{5}}{2}  & \frac{3 - \sqrt{5}}{2}  & \sqrt{5}-2 \emat.
\]

\section{The matrix $A_6$}                       \label{Sinkhorn:section:A_6} 

The construction of the Sinkhorn limit of the $3 \times 3$ matrix
\beq          \label{Sinkhorn:KK}
A_6 = \bmat
K & K & 1\\
K & 1 & 1 \\
1 & 1 & 1 
\emat
\eeq
also requires only high school algebra.  
There exists a unique positive diagonal matrix 
$X  = \diag(x,y,z)$ such that 
\[
S(A_6)  = X A_6 X  =  \bmat
Kx^2 & Kx y & x z  \\
Kx y & y^2 & y z  \\
x z & y z & z^2
\emat 
\]
is a doubly stochastic matrix, and so 
\begin{align}
Kx^2 + Kx y + x z  & =1     \label{Sinkhorn:A6-eqn1}              \\
Kx y + y^2 + y z  & = 1  \label{Sinkhorn:A6-eqn2}             \\
x z + y z + z^2 & = 1.  \label{Sinkhorn:A6-eqn3}            
\end{align}
From~\eqref{Sinkhorn:A6-eqn1} and~\eqref{Sinkhorn:A6-eqn2} we obtain 
\[      
z = \frac{1}{x} -Kx - Ky = \frac{1}{y} -Kx - y 
\]
and so 
\beq                        \label{Sinkhorn:A6-eqn4}  
x = \frac{y}{(K-1)y^2  + 1}
\eeq
and
\beq               \label{Sinkhorn:A6-eqn5}      
z = \frac{1}{y} - \frac{Ky}{(K-1)y^2 + 1} - y
 = \frac{- (K-1)y^4 - 2y^2 + 1}{y( (K-1)y^2+1)}.
\eeq
Inserting~\eqref{Sinkhorn:A6-eqn4} and~\eqref{Sinkhorn:A6-eqn5} into~\eqref{Sinkhorn:A6-eqn3}    
and simplifying, we obtain 
\[
\left( (K-1)y^2 + 1 \right)^3 = K
\]
and so 
\[
y^2 = \frac{K^{1/3} -1}{K-1} = \frac{1}{1+K^{1/3} + K^{2/3}}
\]
and
\[
y = \frac{1}{\sqrt{1+K^{1/3} + K^{2/3}}}.
\]
Inserting this  into~\eqref{Sinkhorn:A6-eqn4} gives 
\[
x = \frac{y}{K^{1/3}} = \frac{1}{ K^{1/3}\sqrt{1+K^{1/3} + K^{2/3}}}.
\]
and then~\eqref{Sinkhorn:A6-eqn5} gives 
\[
z = K^{1/3}y =  \frac{K^{1/3}}{ \sqrt{1+K^{1/3} + K^{2/3}}}.
\]
This determines the scaling matrix X.  
The Sinkhorn limit is the circulant matrix 
\[
S(A_6) = \bmat a & b & c \\ b & c & a \\ c & a & b \emat
\]
with 
\begin {align*}
a  & = Kx^2 = yz = \frac{K^{2/3} -K^{1/3} }{K-1} \\
b & = Kxy = z^2   = \frac{K - K^{2/3} }{K-1} \\
c & = xz = y^2  = \frac{K^{1/3} -1}{K-1}.
\end {align*}
The asymptotic limits are 
\[
\lim_{K\rightarrow \infty} S(A_6) = \bmat 0 & 1 & 0 \\ 1 & 0 & 0 \\ 0 & 0 & 1 \emat 
\qand 
\lim_{K\rightarrow 0} S(A_6) = \bmat 0 & 0 & 1 \\ 0 & 1 & 0 \\ 1 & 0 & 0 \emat. 
\]

\section{The matrix $A_7$}                  \label{Sinkhorn:section:A_7}

Consider the symmetric $3\times 3$ matrix 
\[
A_7= \bmat
K & K & 1 \\
K & 1 & 1 \\
1 & 1 & K
\emat.
\]
There exists a unique positive diagonal matrix $X = \diag(x,y,z)$ such that 
\[
S(A_7) = XA_7X = \bmat
 K x^2 &  K xy & xz \\
 Kxy & y^2 & yz \\
xz & yz & Kz^2
\emat 
\]
is doubly stochastic.  Therefore, 
\begin{align}
 K x^2 +  K xy+ xz  & = 1          \label{Sinkhorn:A7-eqn1}       \\
 Kxy + y^2 + yz  & = 1                 \label{Sinkhorn:A7-eqn2}   \\
xz + yz +  K z^2 & = 1           \label{Sinkhorn:A7-eqn3}  
\end{align}
Because equations~\eqref{Sinkhorn:A7-eqn1} and~\eqref{Sinkhorn:A6-eqn1} 
are identical, and equations~\eqref{Sinkhorn:A7-eqn2} and~\eqref{Sinkhorn:A6-eqn2} 
are identical, we obtain~\eqref{Sinkhorn:A6-eqn4} and~\eqref{Sinkhorn:A6-eqn5}.  
Inserting these formulae for $x$ and $z$ into~\eqref{Sinkhorn:A7-eqn3} gives 
the octic polynomial   
\[
(K-1)^3y^8+3(K-1)^2y^6-(K-1)(2K-3)y^4-(4K-1)y^2+K = 0.
\]
By Theorem~\ref{Sinkhorn:theorem:Sinkhorn-1}, this polynomial has at least one  solution $y \in (0,1)$.
If $K>1$, then, by Descartes's rule of signs, this polynomial has exactly two positive solutions.  
If $0 < K < 1$, then this polynomial has one or three positive solutions.  
For matrices of the form $A_7$, we do not have explicit formulae for the coordinates 
of the Sinkhorn limit as functions of $K$.  Computer calculations suggest 
that the asymptotic limits of $S(A_7)$ as $K\rightarrow \infty$ and $K\rightarrow 0 $ are 
\[
 \bmat 0 & 1 & 0 \\ 1 & 0 & 0 \\ 0 & 0 & 1 \emat 
 \qand
  \bmat 0 & 0 & 1 \\ 0 & 1 & 0 \\ 1 & 0 & 0 \emat.  
\]

\section{Gr\" obner bases and algebraic numbers}              \label{Sinkhorn:section:Grobner}
I like solving problems using high school algebra.  
However, it is important to note that the previous calculations are also easily 
done using Gr\" obner bases.

For every $n \times n$ matrix $A = (a_{i,j})$ and diagonal matrix $X = \diag(x_1,\ldots, x_n)$, 
we have the matrix 
\[
XAX = \bmat a_{i,j} x_i x_j \emat.
\]
If $A$ is positive and symmetric, then, by Theorems~\ref{Sinkhorn:theorem:Sinkhorn-1} 
and~\ref{Sinkhorn:theorem:Sinkhorn-2}, the $n$ quadratic equations  
\[
q_i = q_i(x_1,\ldots, x_n) = \sum_{j=1}^n a_{i,j} x_ix_j - 1 = 0 \qquad \qquad (i = 1,\ldots, n) 
\]
have a unique positive solution, and the diagonal matrix 
$X = \diag(x_1,\ldots, x_n)$ is the unique scaling matrix in the Sinkhorn limit  $S(A) = XAX$.
Equivalently, $(x_1,\ldots, x_n)$ is the unique positive vector in the affine variety 
of the ideal in $\R[x_1,\ldots, x_n]$ generated by the set of polynomials $\{q_1,\ldots, q_n\}$.  
For each lexicographical ordering of the variables $x_1,\ldots, x_n$,
Maple (and other computer algebra programs) can compute a Gr\" obner basis for the ideal.  
The Gr\" obner basis for this ideal shows that  if the coordinates of the matrix $A = (a_{i,j})$ 
are rational numbers, then $x_1,\ldots, x_n$ are algebraic numbers of degrees 
bounded in terms of $n$.

Here is an example.   Let $n = 3$ and $X = \diag(x,y,z)$. 
Consider the matrices 
\[
A_7 = \bmat
K & K & 1 \\
K & 1 & 1 \\
1 & 1 & K
\emat
\qqand
XA_7X = 
\bmat
K x^2& Kxy & xz \\
Kxy & y^2 & yz \\
xz & xy & Kz^2
\emat.
\]
with $K>0$ and $K \neq 1$.  
There exist unique positive real numbers $x,y,z$ that satisfy the quadratic  
equations 
\begin {align*}
 K  x^2 +  K  xy+ xz  & = 1\\
 K xy + y^2 + yz & = 1 \\
xz + yz +  K  z^2 & = 1. 
\end {align*}
Equivalently, $(x,y,z)$ is the unique positive vector in the affine variety 
$V(I)$, where $I$ is the ideal in $\R[x,y,z]$ generated by the polynomials 
\begin {align*}
 K  x^2 +  K  xy+ xz  & - 1\\
 K xy + y^2 + yz  & - 1 \\
xz + yz +  K  z^2 & - 1. 
\end {align*}
Let $K=2$.  
Using the Groebner package in Maple with the lexicographical order $(x,y,z)$, we obtain 
the Gr\" obner basis 
\begin {align*}
f_1(z) & = 4-28 z^2+62 z^4-57 z^6+18 z^8 \\
f_2(y,z) & =  -17 z^3+39 z^5-18 z^7+2 y \\
f_3(x,z) & =   -20 z+96 z^3-135 z^5+54 z^7+4 x.
\end {align*}  
Applying Maple with the lexicographical order $(y,z,x)$, we obtain 
the Gr\" obner basis 
\begin {align*}
g_1(x) & = 2-17 x^2+22 x^4+48 x^6+36 x^8 \\
g_2(x,z) & = -103 x+378 x^3+624 x^5+396 x^7+14 z \\
g_3(x,y) & =  3 x-56 x^3-72 x^5-36 x^7+7 y.
\end {align*}  
Applying Maple with the lexicographical order $(z,x,y)$, we obtain 
the Gr\" obner basis 
\begin {align*}
h_1(y) & = 2-7 y^2-y^4+3 y^6+y^8 \\
h_2(x,y) & =  -4 y+2 y^5-3 y^3+y^7+6 x \\
h_3(y,z) & =   -7 y+5 y^5+3 y^3+y^7+6 z. 
\end {align*}  
Thus, $x^2$, $y^2$, and $z^2$ are algebraic numbers of degree at most 4, and we have explicit polynomial representations of each variable $x$, $y$, $z$ in terms of the other two variables.

 For arbitrary $K$, applying Maple with the lexicographical order $(y,z,x)$, we obtain 
the Gr\" obner basis 
\begin {align*}
h_1(y) 
& = K- (4 K - 1) y^2 -(K-1)(2K-3) y^4+3 (K-1)^2  y^6+(K-1)^3 y^8 \\
h_2(x,y) 
& =   K (K+1) x  -2 K y - (K-1)(2K-1) y^3+ 2 (K-1)^2  y^5+(K-1)^3 y^7 \\
h_3(y,z) 
& = K (K+1) z - (K-1)^2  y -3(K-1) y^3+ (K-1)^2(K-3) y^5+( K-1)^3 y^7.
\end{align*}   
For each of the 8 roots of $h_1(y)$, the polynomials $h_2(z,y)$ and $h_3(x,y)$ 
determine unique numbers $x$ and $z$.  
Exactly one of the triples $(x,y,z)$ will be positive.


\section{Rationality and finite length}      \label{Sinkhorn:section:finiteLength}
For what positive $n\times n$ matrices does the alternate scaling algorithm 
converge in finitely many steps?   This problem has been solved for $2 \times 2$ matrices 
(Nathanson~\cite{nath2019-184}), but it is open for all dimensions $n \geq 3$.  
In dimension 3, matrices equivalent to $A_1$ become doubly stochastic in one step, 
that is, after one row or one column scaling.  
Ekhad and Zeilberger~\cite{ekha-zeil19} computed a positive  $3\times 3$ matrix
 that becomes doubly stochastic in exactly two steps, and Nathanson~\cite{nath2019-186} 
generalized this construction.  It is not know if there exists  a positive  $3\times 3$ matrix 
that becomes doubly stochastic in exactly $s$ steps for some $s \geq 3$.

Consider the matrix $A_2 = \bmat K & 1 & 1 \\ 1 & 1 & 1 \\ 1 & 1 & 1 \emat$ with parameter $K$.   
If $K$ is a rational number, then every matrix generated by iterated row 
and column scalings has rational coordinates.  
If the Sinkhorn limit contains an irrational coordinate, 
then the alternate scaling algorithm cannot terminate in finitely many steps.  

Let $K$ be an integer,  $K \geq 2$.
In Section~\ref{Sinkhorn:section:A_2}  we proved that the Sinkhorn limit $S(A_2)$ 
has coordinates in the quadratic  field 
 $\Q(\sqrt{8K+1})$.  For example, from~\eqref{Sinkhorn:A2-a},  
 the $(1,1)$ coordinate of $S(A_2)$ is 
 \[
 \frac{2K+1-\sqrt{8K+1}}{2(K-1)}.
 \]
This number is rational if and only if the odd integer $8K+1$ is the  square of an odd integer, that is,
if and only if $8K+1 = (2r+1)^2$ for some positive integer $r$ 
and  so $K = r(r+1)/2$ is a triangular number.  
From~\eqref{Sinkhorn:A2-a},~\eqref{Sinkhorn:A2-b}, and~\eqref{Sinkhorn:A2-c}, 
we obtain 
\begin {align*}
a & = \frac{r^2-r}{r^2+r-2} = \frac{r}{r+2}  \\
b & = \frac{r -1}{r^2+r-2}  = \frac{1}{r+2} \\
c & =   \frac{r^2-1}{2(r^2+r-2)} = \frac{r+1}{2(r+2) }.
\end {align*}
Moreover, $S(A_2) = X A_2 X$, where $X = \diag (x,y,y)$ 
with $Kx^2 = a$ and $y^2 = c$.  Thus, 
\[
x = \sqrt{\frac{a}{K} } =  \sqrt{ \frac{2}{ (r+1)(r+2)}  }
\qqand 
y = \sqrt{c} =  \sqrt{   \frac{r+1}{2(r+2)}    }.  
\]

For example, if $K = 3$, then $r=2$ and 
\[
A_2 = \bmat 3 & 1 & 1 \\ 1 & 1 & 1  \\ 1 & 1 & 1 \emat
\qqand  
X A_2 X = S(A_2 ) = \bmat 1/2 & 1/4 & 1/4 \\ 1/4 & 3/8 & 3/8  \\ 1/4 & 3/8 & 3/8 \emat 
\]
where 
\[
 X = \diag( \sqrt{6}/6, \sqrt{6}/4,  \sqrt{6}/4 ).
\]
Note that  $A_2$ also has a scaling by rational matrices  
\[
S(A_2) = X'A_2 Y'
\]
where 
\[
X' = \diag( 1/6, 1/4,1/4) \qqand Y' = \diag(1, 3/2, 3/2).
\]
It is not known if there exists a triangular number $K$ 
for which the alternate scaling  algorithm terminates in a finite number of steps.

\section{Open problems}

\benum
\item
Compute explicit formulas for the Sinkhorn limits of  matrices of the form $A_7$.
More generally, compute explicit formulas for the Sinkhorn limits of all positive symmetric 
$3\times 3$ matrices.
This is a central problem.

\item
Here is a special case.  
Let $K, L, M$ and 1 be pairwise distinct positive numbers.  Compute the Sinkhorn limits of the matrices
\[
\bmat 
K & 1 & 1 \\
1 & L & 1 \\
1 & 1 & 1
\emat
\qqand
\bmat 
K & 1 & 1 \\
1 & L & 1 \\
1 & 1 & M
\emat.
\]

\item
For what positive $n\times n$ matrices does the alternate scaling algorithm 
converge in finitely many steps?   This is the problem discussed in Section~\ref{Sinkhorn:section:finiteLength}.

\item
It is not known what algebraic numbers appear as coordinates of Sinkhorn limits 
of matrices with positive integral coordinates.  
It would be interesting to have an example of an algebraic number in the unit interval 
that is not a coordinate of the Sinkhorn limit of a positive integral matrix.

\item
Does every possible shape  of a doubly stochastic $3\times 3$ matrix $A$ appear as the 
nontrivial Sinkhorn  limit of some $3\times  3$ matrix?

\item
Why does the shape of the Sinkhorn limit $S(A)$ seem to depend only on the shape of the matrix $A$ 
and not on the numerical values of the coordinates of $A$?

\item
Let A\ be a nonnegative $m\times n$ matrix.
Let $ \mbr = (r_1, r_2,\ldots, r_m) \in \R^m$ and let $ \mbc = (c_1, c_2, \ldots, c_n) \in \R^n$.
The matrix A\ is \emph{$ \mbr$-row stochastic} 
if $\rowsum_i( A ) = r_i$ for all $i \in \{1,2,\ldots, m\}$.
The matrix A\ is \emph{$ \mbc$-column stochastic} 
if $\colsum_j( A ) =c_j$ for all $j \in \{1,2,\ldots, n\}$.
The matrix $ A $ is \emph{$(\mbr,\mbc)$-stochastic}  
if it is both $\mbr$-row stochastic and $\mbc$-column stochastic.

Let A\ be a positive matrix.  Let $X$ be the $m\times m$ diagonal matrix whose 
$i$th coordinate is $r_i/\rowsum_i(A)$, and let $Y$ be the $n \times n$ diagonal matrix whose 
$j$th coordinate is $c_j/\colsum_j(A)$.
The matrix $XA$ is $\mbr$-row stochastic and the matrix 
$AY$ is $\mbc$-column stochastic.  
A simple modification of the alternate scaling algorithm produces 
an $(\mbr,\mbc)$-stochastic Sinkhorn limit.  
It is an open problem to compute explicit Sinkhorn limits in 
the $(\mbr,\mbc)$-stochastic setting.

\item
It is a old problem in number theory to understand the continued fractions 
of the cube roots of integers, and, in particular, to understand 
the approximation of $\sqrt[3]{2}$ by rationals.  
One coordinate of the Sinkhorn limit of the matrix $A_6$  with $K=2$  is $\sqrt[3]{2} - 1$.  
The matrix $A_6$  with $K=2$  has rational coordinates, and so 
the matrices constructed by the alternate scaling algorithm also have 
rational coordinates, and generate explicit sequences of rational 
approximations to $\sqrt[3]{2}$.   The nature of these approximations
remains mysterious. 
\eenum

\section{Notes}
The computational complexity of Sinkhorn's alternate scaling algorithm 
is investigated in 
Kalantari and Khachiyan~\cite{kala-khac93,kala-khac96},
Kalantari,  Lari, Ricca, and Simeone~\cite{kala08}, 
Linial, Samorodnitsky and Wigderson~\cite{lini-samo-wigd98b}
and Allen-Zhu, Li, Oliveira, and Wigderson~\cite{alle-li-oliv-wigd17}.
An extension of matrix scaling to operator scaling began 
with Gurvits~\cite{gurv04}, and is developed in 
Garg,  Gurvits, Oliveira, and Wigderson~\cite{garg-gurv-oliv-wigd16,garg-gurv-oliv-wigd17},
Gurvits~\cite{gurv15}, and Gurvits and Samorodnitsky~\cite{gurv-samo14}.  
Motivating some of this recent work are the classical papers of Edmonds~\cite{edmo67} 
and Valient~\cite{vali79a,vali79b}.

The literature on matrix scaling is vast.  See the recent survey paper of 
Idel~\cite{idel16}.  For the early history of matrix scaling, see  
Allen-Zhu, Li, Oliveira, and Wigderson~\cite[Section 1.1]{alle-li-oliv-wigd17}.

\emph{Acknowledgements.}
The alternate scaling algorithm was discussed in several lectures 
in the New York Number Theory Seminar, 
and I thank the participants for their useful remarks.  
In particular, I thank David Newman for making the initial computations 
that suggested some of the problems considered in this paper.  
I also benefitted from a careful and thoughtful referee's report.  

\def\cprime{$'$} \def\cprime{$'$} \def\cprime{$'$}
\providecommand{\bysame}{\leavevmode\hbox to3em{\hrulefill}\thinspace}
\providecommand{\MR}{\relax\ifhmode\unskip\space\fi MR }
\providecommand{\MRhref}[2]{%
  \href{http://www.ams.org/mathscinet-getitem?mr=#1}{#2}
}
\providecommand{\href}[2]{#2}

\end{document}